\documentclass[11pt,a4paper]{article}

\usepackage[margin=2.5cm]{geometry}
\usepackage{amsmath,amssymb}
\usepackage{booktabs}
\usepackage{array}
\usepackage{enumitem}

\usepackage{titlesec}
\titleformat{\section}{\large\bfseries}{\thesection}{1em}{}
\titlespacing*{\section}{0pt}{1.5ex plus 0.5ex minus 0.2ex}{1ex plus 0.2ex}
\titleformat{\subsection}{\normalsize\bfseries}{\thesubsection}{1em}{}
\titlespacing*{\subsection}{0pt}{1.2ex plus 0.3ex minus 0.2ex}{0.8ex plus 0.2ex}

\usepackage{hyperref}
\hypersetup{
  colorlinks=true,
  linkcolor=blue,
  citecolor=blue,
  urlcolor=blue
}

\title{\textbf{\Large The Four-Element Theory of Structure-Driven Methodology: \\ A Methodological Shift from the Perspective of Natural Science and Engineering Practice}}
\author{
    {CHEN Shengchang} \\
    School of Earth Sciences, Zhejiang University, Hangzhou, China \\
    \texttt{chenshengc@zju.edu.cn}
}
\date{}

\begin{document}

\maketitle

\begin{abstract}
This paper directly addresses the systemic structural blind spot caused by ``objective pre-setting'' in contemporary scientific problem-solving, and extracts an emerging methodological shift from concrete practices in natural sciences and engineering. Solving activities have long been fundamentally Objective-Driven: an evaluative objective is first set based on experience, and then approximated through iterative trial and error. While this path is reliable for simple linear problems, it easily falls into local optima, yields non-reproducible results, and distorts native attributes when dealing with complex nonlinear problems. This paper proposes the Four-Element Theory of Structure-Driven Methodology, reconstructing the solving process into four sequentially irreversible stages: ``Deconstruction -- Structure -- Construction -- Solving''. The Deconstruction Element suspends \textit{a priori} presuppositions, reducing the problem to mutually independent basic units. The Structure Element penetrates superficial information to extract the rigid invariants inherent to the problem that remain unchanged regardless of observation or intervention. The Construction Element uses these invariants as anchors to construct an equivalent mapping from problem inputs to the solution. The Solving Element executes the mapping to directly generate the final solution. The fundamental basis for this four-element division is not logical symmetry, but the distinction of purpose: deconstruction, structure, and construction all prepare for solving, while the solving element is the actual implementation stage of problem resolution. Evaluation metrics and experience are relegated to the position of \textit{a posteriori} verification. Based on empirical foundations from various practices---such as number theory proposition proofs, long-term evolution simulation of Hamiltonian systems, deep learning signal processing and structure-driven backpropagation, and seismic wave inversion imaging---the paper demonstrates the cross-domain consistency of this methodology. Structure-Driven methodology does not negate the Objective-Driven approach, but complements it: when structure is sufficient, solutions are directly constructed and solved; when structure is lacking, it degenerates into the Objective-Driven path. The specific form of the Four-Element Theory in the field of inverse problems is Structure-Driven Inversion (SDI), which includes Mathematical-Structure-Driven Inversion (MSDI) and Physical-Structure-Driven Inversion (PSDI). Spectral-domain unitary diagonalized pseudo-inverse, true 3D gravity inversion, Physical-Structure-Driven Waveform Inversion Imaging (PSDWII), and Physical-Structure-Driven Backpropagation (PSDBP) are all successful practices of the Four-Element Theory in different scenarios. The validity of this methodology depends on identifiable and verifiable rigid invariants; it degenerates to the Objective-Driven approach in the absence of structure.

\vspace{1em}
\end{abstract}

\vspace{1em}


\section{Introduction: A Familiar Dilemma}

We have all had this experience: when facing a complex problem, we first set a ``good'' standard based on experience---error should be small, cost should be low, effect should be optimal---and then continuously adjust the path step by step to approach this standard. If the problem is simple, this method is very effective. But if the problem is complex, troubles arise: no matter how much you adjust, it always falls just short; if someone else does it, the result is different; sometimes, in order to meet the standard, the problem itself gets distorted.

This dilemma does not only appear in daily life. In science and engineering, it has a more formal name: ``Objective-Driven'' or ``Optimization-Driven''. Its basic logic is: first set a comparable evaluative objective based on experience, and then through trial and error, adjustment, and iteration, make the candidate solution gradually approach the objective. From Laplace's least absolute deviation criterion, Gauss and Legendre's least squares method, to Tikhonov regularization \cite{ref2}, the adjoint state method \cite{ref3,ref4}, and contemporary large-scale optimization algorithms, the form of objective functions and optimization algorithms have changed over the past two hundred years, but the core logic of ``constructing an objective function and optimizing it'' has not.

The advantage of Objective-Driven lies in its simplicity. As long as the computational process can be carried out, any problem---whether differentiable or not---can be incorporated into the framework of objective function optimization. However, the cost of simplicity lies in simplicity itself: a simple framework capable of accommodating all problems inevitably cannot be tailored to specific ones. When a process with translational invariance is treated as a general mapping, when a physical process with a causal chain is treated as a general fitting problem, the problem's own structure is obscured.

This leads to three limitations. First, the efficiency bottleneck: every iteration requires a complete calculation, and the computational cost for high-dimensional complex problems expands dramatically. Second, blurred meaning: changes in metrics are merely mathematical constructs; they indicate ``how parameters should change to lower the objective,'' but do not directly explain the causal relationship between parameters and observations. Third, coupling and distortion: in multi-parameter problems, the contributions of different factors are easily mixed, creating artificial coupling; to further reduce metrics, the solving process may introduce false constraints inconsistent with the nature of the problem.

These limitations are not imperfections in engineering implementation, but logical consequences of the Objective-Driven methodology: as long as the solving path is anchored to externally preset objectives rather than the problem's own structure, they cannot be eliminated at the root.

In recent years, a different methodology has been independently emerging in multiple fields. In 2026, a formal proof based on automated reasoning of large models was released on GitHub, proving a weakened form of the Goldbach conjecture for the Liouville function; the proof was only two pages, and the core was not exhaustive search, but identifying the sign structure of the Liouville function under additive combinations, and directly constructing the proof path along the structure. Wang Hong and Zahl's proof of the three-dimensional Kakeya conjecture \cite{ref5} also first deconstructed the original problem into ``the degree of overlap of thin tubes,'' identified structures such as multi-scale nesting and convex set density, and then constructed an inductive path along the structure. Noether's theorem in physics \cite{ref6} revealed long ago: symmetry dictates conservation laws, structure precedes solving. In the field of pure mathematics, MIT mathematician Larry Guth similarly points out that the essence of mathematical research lies in discovering deep connections and structures between seemingly unrelated things (``Things come together''). Especially in the face of the computational impact of AI, human mathematicians need to define the true meaning of their research through reflection \cite{ref1}.

In the fields of deep learning and geophysics, several core research works by the author in recent years have repeatedly pointed to the same logic: in inverse problem solving, unitary transformation is used to deconstruct the problem, identifying its quasi-Singular Value Decomposition (quasi-SVD) structure, and then constructing a spectral-domain pseudo-inverse \cite{ref7} to directly solve the true 3D gravity inversion; in seismic wave inversion imaging, Physical-Structure-Driven Waveform Inversion Imaging (PSDWII) was proposed \cite{ref8}, utilizing virtual source projection to directly construct the imaging mapping; in neural network training, the Physical-Structure-Driven Backpropagation (PSDBP) framework was proposed \cite{ref9}, reconstructing gradient search into structure propagation; in signal denoising, the Signal-Fitting Diffusion Denoising Probabilistic Model (SFDDPM) was led and developed \cite{ref10}. These practices span mathematical inversion, physical wavefields, and deep learning, yet all point to the same logic: identify structure, utilize structure, and construct the solving path.

The task of this paper is to extract a common methodological structure from these concrete experiences in natural science and engineering practice. It is neither a purely philosophical speculative essay nor a purely technical report, but rather abstracts the same set of execution logic---Deconstruction, Structure, Construction, Solving---from different practices such as mathematics, physics, deep learning, and geophysics. This paper calls it the ``Four-Element Theory of Structure-Driven Methodology.''

It is necessary to clarify the relationship between the ``Four-Element Theory'' proposed in this paper and the author's previous ``Structure-Driven Inversion (SDI)'' work. SDI is a new paradigm proposed for inverse problems \cite{ref15}, containing two paths: Mathematical-Structure-Driven Inversion (MSDI) and Physical-Structure-Driven Inversion (PSDI). The Four-Element Theory proposed in this paper generalizes the logic of SDI to a macro-level meta-methodological framework for general problem-solving; \textbf{SDI is incorporated into it, becoming a specific form of the Four-Element Theory in the field of inverse problems.} It should be particularly emphasized that whether it is PSDWII based on seismic wave propagation (a true physical process) or PSDBP based on forward propagation of neural network information (a quasi-physical process), their essence is to deconstruct the propagation process, identify its causal structure, and then construct a propagation mapping without blind search. They both belong to the successful practices of the Four-Element Theory in physical and quasi-physical problems.

This paper first analyzes the logical chain triggered by ``objective pre-setting'' in the Objective-Driven methodology; then proposes the four stages of the Four-Element Theory, defining the input, output, and source of legitimacy for each stage; demonstrates the application of the framework through cross-disciplinary cases; compares it with existing methodologies; discusses the applicable boundaries and risks of this framework; and finally provides conclusions and future research directions.

\section{Why Does the Objective-Driven Paradigm Fall into Dilemma?}

To understand the necessity of Structure-Driven methodology, one must first clarify the causal chain of the Objective-Driven paradigm. The starting point of this chain is not iteration itself, but the pre-setting of the evaluative objective prior to structure identification.

\subsection{Objective in the Forefront, Structure in the Background}
The most fundamental feature of Objective-Driven is the pre-setting of the evaluative objective. Before solving begins, the solver has already determined the standard for what constitutes a good solution based on experience. This pre-setting constraint has two forms.

The first is the implicit form. In daily problem-solving, the evaluative objective may not be explicitly written out, but acts in the form of common sense, convention, or empirical judgment. The solver has a ruler in mind, but this ruler has never been examined.

The second is the explicit form. In scientific computing and engineering optimization, the evaluative objective is formalized into objective functions, performance metrics, loss functions, or evaluation scales. Formalization brings computability, but also stronger constraints: once the objective function is written down, the solving path is predetermined by it---gradient directions, regularization forms, and convergence criteria are all determined by the objective function.

The common point of both forms is: the evaluative objective temporally precedes structure identification, is logically independent of it, and even replaces it. Before the solver has time to ask ``what structure does this problem itself have,'' they have already been drawn by the objective into the cycle of ``how to get closer to the objective.''

\subsection{Structural Blind Spot: The Objective is Like a Beam of Light, Illuminating Only a Part}
The reason pre-set objectives constitute a problem is their systemic obscuring effect on the internal structure of the problem.

The so-called internal structure refers to the network of relationships inherent to the problem itself, which does not change with the solver's will or observational data. It may be symmetries at the operator level, orthogonality at the mapping level, conservation laws at the system level, or irreducible constraint relationships between variables. Such structures are the essential attributes of the problem, independent of any evaluation system.

The evaluative objective, however, is an external construct set by the solver based on experience. When solving is anchored solely to the evaluative objective, the problem is projected onto the coordinate axes of the evaluation standard: factors that reduce metrics are retained, while factors unrelated to metrics are ignored. Structural information carried by the problem itself, if it cannot be translated into gains in metrics, will be filtered out. This is the structural blind spot.

To use an analogy: the objective is like a beam of light, illuminating only the parts related to metrics, while the problem's own structure falls into darkness. This is not because the solver's eyesight is poor, but because it is determined by the direction of the light source.

The structural blind spot is not an oversight by the solver, but an inherent logical consequence of a methodology anchored to external objectives. Nor does it mean the solver cannot see the structure. On the contrary, domain experts often have profound intuition about problem structures. The problem is that the Objective-Driven methodology provides no methodological position for this intuition: it can only influence regularization terms, initial models, or search strategies in the name of experience, but cannot become the construction basis for the solving path. Structure is downgraded to a part of experience, and experience can only enter the solving process in the form of constraints, not construction.

\subsection{Iterative Trial and Error Becomes the Only Way Out}
Once the structural blind spot is formed, iterative trial and error becomes the only feasible path.

Since the solving path cannot be directly given by the problem's structure, it can only be gradually approached through cycles of trial, feedback, and correction. In this cycle, the evaluative objective acts as the feedback signal, gradients or heuristic rules act as the correction direction, and the number of iterations acts as the cost of approximation. The solving process thus degenerates into a search activity of repeatedly finding paths in unknown terrain.

This brings threefold costs. First, path locking and local optima: the search starting point and step size are set by experience, and once trapped in a local extremum, the gradient signal disappears and the solving process is locked. Second, non-reproducible results: differences in experience among different solvers are reflected in initial models, parameters, and stopping criteria, making it difficult to stably reproduce results for the same problem. Third, distorted native attributes: to further reduce metrics, the solving process may introduce false constraints inconsistent with the nature of the problem.

Thus, the complete causal chain of the Objective-Driven methodology is clarified: experience-constrained objective evaluation leads to structural blind spots, structural blind spots lead to iterative trial and error, and iterative trial and error leads to local optima, non-reproducibility, and attribute distortion. The starting point of the chain is not iteration, but the pre-setting of the evaluative objective; to break this chain, the anchor must be changed.

\subsection{A Clarification: Structural Blind Spot Does Not Equal Lack of Knowledge}
It should be emphasized that the \textbf{structural blind spot} proposed in this paper belongs to an information filtering effect at the methodological level, rather than referring to a lack of domain knowledge among individual researchers. In many engineering and scientific practices, even if researchers have mastered rich professional knowledge and clearly understand the mathematical constraints and physical laws behind the problem, they may still fall into the structural blind spot. This phenomenon is universal: experts can qualitatively identify the symmetry, conservation laws, or inherent properties of operators possessed by the system, but when carrying out actual solving work, the entire process is still constrained by the Objective-Driven framework. At this point, the acquired structural knowledge cannot be elevated to the generative basis of the solving path, and can only be downgraded to auxiliary empirical means, such as being transformed into regularization terms, prior constraints, initial models, or parameter tuning experience. Although such knowledge can constrain the solution space and suppress the ill-posedness of the problem, it cannot directly determine the solving mapping itself.

In other words, researchers may already possess relevant knowledge, but the solving methodology limits the position where knowledge plays its role. Under the Objective-Driven methodology, external evaluation objectives occupy the core anchor of solving; even if researchers master correct structural knowledge, this knowledge can only participate in solving as a ``soft constraint,'' and cannot serve as a rigid basis for constructing solutions. Expanding knowledge reserves certainly helps to set more reasonable regularization conditions and optimize hyperparameters, improving the quality of solutions within a certain range, but it cannot fundamentally cut off the causal chain of ``objective pre-setting $\rightarrow$ structural blind spot $\rightarrow$ iterative trial and error.''

Therefore, the way to overcome the structural blind spot is not to accumulate more domain knowledge, but to change the anchor of the solving path. Domain knowledge can improve the quality of the solution, but only moving the anchor can change the solving logic.

\section{The Four-Element Theory of Structure-Driven Methodology: Four Stages}

The Four-Element Theory reconstructs the solving process into four sequentially progressive and irreversible stages: Deconstruction, Structure, Construction, Solving. The four respectively respond to four questions: What is the true composition of the problem? What is the rigid structure inherent to the problem? How to construct a solving mapping based on the structure? How to execute the mapping and generate the final solution?

The basis for the four-element division is not logical symmetry, but the distinction of purpose: deconstruction, structure, and construction all prepare for solving, and the solving element is the implementation stage of problem resolution. Separating construction and solving is to prevent solving from being replaced by objective functions or evaluation metrics, ensuring that solving directly points to the problem itself.

\subsection{Deconstruction Element: Reducing the Problem without Prejudgment}
The Deconstruction Element is the starting point of solving. Its task is not simple decomposition, but a methodological suspension: suspending all \textit{a priori} assumptions, preset objectives, and evaluation standards from past experience, and facing the complete problem without disciplinary bias or empirical prejudgment.

The goal of deconstruction is to reduce the problem to several basic units, and these units must satisfy strict independence requirements: no cross-coupling terms exist between them, no hidden linkage relationships exist, and a state change in any one unit does not cause a concomitant change in the states of other units. This requirement is far stronger than decomposition in the usual sense. General decomposition allows information exchange between sub-problems, while the Deconstruction Element requires the complete stripping of coupling.

Why is it so strict? Because subsequent structure identification must occur between units. If coupling remains between units, the identified relationships will be mixed with pseudo-relationships caused by artificial decomposition, and the rigidity of the structure cannot be guaranteed. The strictness of deconstruction directly determines the reliability of the Structure Element.

The Deconstruction Element does not preset a solving objective. This non-prejudgment requirement may seem passive, but is actually the logical premise of the entire methodology: only on a complete, authentic problem that has not been projected by evaluation standards can the internal structure be fully identified. Once objective pre-setting is mixed into the deconstruction stage, the structural blind spot will reappear at the starting point.

In practice, deconstruction presents two basic types. The first is spatial deconstruction, which divides the problem into several units according to the natural spatial or object boundaries of the problem, such as discretizing a continuous medium into grid nodes, splitting an entire system into functional modules, or splitting a composite signal into single sources. The second is frequency-domain deconstruction, which splits the problem into independent components according to frequency or scale components, such as decomposing a wavefield into plane waves of different wavenumbers, or decomposing a signal into harmonics of different frequencies. The common requirement of both types of deconstruction is consistent: the units after splitting must be mutually independent.

It is also necessary to distinguish the Deconstruction Element from so-called problem simplification. Problem simplification aims to reduce solving difficulty, using methods such as ignoring secondary factors, reducing dimensions, and adopting approximate models, at the cost of introducing model errors; the Deconstruction Element aims to strip coupling, does not ignore any factors, does not reduce dimensions, but only changes the combination of factors to separate originally mixed relationships. Therefore, deconstruction does not introduce approximation errors, which is the fundamental reason why it can serve as the premise for structure identification.

\subsection{Structure Element: Finding That Which Never Changes}
The Structure Element is the core of the Four-Element Theory, and the key that distinguishes it from all experience-driven paths.

The task of the Structure Element is to penetrate superficial information on the basic units obtained from deconstruction, and extract the rigid invariants inherent to the problem itself. The so-called rigidity has three meanings: it does not change with the values of observational data, does not change with the solver's subjective intervention, and does not change with the specific scenario in which the problem is situated. It belongs to the problem itself, not to a particular instance of solving.

The criterion of rigidity lies in whether it remains invariant under the above three types of perturbations. In physics, for example, the symplectic structure of a Hamiltonian system remains invariant under any canonical transformation; in mathematics, the spectral decomposition of a unitarily diagonalizable operator remains invariant under any unitary transformation; in number theory, the complete multiplicativity of the Liouville function remains invariant under any integer factorization; in deep learning, the causal chain of network forward propagation has a decomposable structural invariance under parameter updates. Such invariances are not statistical regularities derived from data, but constraints directly derived from the problem definition.

In terms of types, rigid invariants can be divided into two major categories. The first is mathematical structure, such as the unitarity of an operator, the orthogonality of a mapping, the diagonalizability of a matrix, and the commutative and associative laws of algebraic operations. The second is physical structure, such as the symmetry of a system, conservation laws, the decomposability of propagation equations, and physical constraints between parameters. Mathematical structure determines whether direct solving is possible, and physical structure determines whether the solution has physical meaning.

In practice, rigid invariants may appear as explicit analytical constraints (such as known wave equations or conservation laws), or as implicit algebraic or topological structures (such as group symmetries or spectral properties of operators). Explicit structures are relatively easy to extract, while implicit structures often require deep mathematical tools (such as functional analysis or topology) to be made explicit, which is the main challenge faced by structure identification in practice.

In the context of inverse problems, as a specific form of the Four-Element Theory, SDI's internal MSDI path utilizes the unitary diagonalizability of operators to directly construct pseudo-inverses in the spectral domain; its PSDI path utilizes the causal chain of wave propagation to solve directly through virtual source projection and deconvolution. In the context of deep learning, SFDDPM utilizes the low-dimensional structure of the signal itself to change the fitting objective, and PSDBP utilizes the structural invariants of network forward propagation to reconstruct the backpropagation path.

The distinction between the Structure Element and empirical common sense must be clarified. Empirical common sense is inductive knowledge about how things are usually done, which changes with domain, individual, and time; rigid invariants are deductive constraints about how the problem necessarily is, which do not change with any of the above factors. The former can serve as inspiration, the latter can serve as the basis for construction. The Objective-Driven methodology misplaces the latter as the former, which is the epistemological root of the structural blind spot.

Although structure identification depends to a considerable extent on domain experience, it is not entirely based on intuition, and several operational paths can be summarized. The first is symmetry analysis, which examines under what transformations the problem remains invariant. The second is spectral analysis, which examines the spectral properties of operators or mappings. The third is invariant testing, which examines quantities that are always conserved under various transformations. The fourth is coupling analysis, which examines whether the coupling between variables can be completely eliminated through some transformation. The four paths are not mutually exclusive, and their common direction is: to separate the structure from the superficial appearance and express it in a verifiable form.

\subsection{Construction Element: Designing a Direct Path from Problem to Solution}
The Construction Element is the implementation and design stage of the Four-Element Theory. Its task is to construct an equivalent mapping, inverse operator, or solving scheme from the problem input to the final solution, based on the extracted rigid invariants. The so-called construction means that the form of the solving path is uniquely determined by the structure, rather than by search. The Construction Element does not execute the final calculation, nor does it guide construction with evaluation metrics; it is only responsible for translating the structure into an executable mapping.

Directness is the essential feature of the Construction Element. The so-called directness means that each step of the solving mapping is determined by the structure, rather than by search. Taking linear inverse problems as an example, when the forward operator has a unitarily diagonalizable structure, its pseudo-inverse can be directly constructed in the spectral domain through unitary diagonalization, without any iteration. At this time, the core operations constructed are division, inversion, or deconvolution, while the core operations of the Objective-Driven path are cross-correlation and multiplication. The difference between the two is not a difference in efficiency magnitude, but a difference in logical type: division naturally separates the contributions of each component, so no mathematical coupling is generated; cross-correlation mixes the contributions of different components together, and mathematical coupling arises from this.

The Construction Element does not absolutely exclude iteration. In cases where the structure is incomplete, stepwise structural projection may still be necessary. But at this time, the nature of iteration has fundamentally changed: each step is a projection with a structural basis, not a blind search. Iteration changes from a means of search to a sequence of structural projections, which is the essence of the methodological shift.

Based on this, two forms of the Construction Element can be distinguished. The first is analytical construction: when the structure is completely sufficient, a closed-form mapping is constructed, and the solution can be directly given by the structure in closed form. The second is projection construction: when the structure is partially sufficient, a sequence of several structural projections is constructed, and the projection direction of each step is determined by the structure. The difference between the two forms lies in the number of steps, not in nature: whether one step or multiple steps, each step has a clear structural basis.

\subsection{Solving Element: Actually Computing the Solution}
The Solving Element is the purpose stage of the Four-Element Theory. Its task is to execute the equivalent mapping given by the Construction Element, completing point-by-point operations or stepwise projections in the transform domain, spectral domain, or projection sequence, directly generating the final solution. The Solving Element does not perform search, does not adjust the objective function, nor does it redefine the path based on experience. Its legitimacy comes entirely from the Construction Element: the Construction Element obtains legitimacy from the structure, and the Solving Element obtains legitimacy from the construction.

The Solving Element is the direct implementation stage of the purpose of ``solving the problem'' in the Four-Element Theory. Deconstruction, structure, and construction all prepare for solving; if solving is not independent, construction easily swallows solving again, and solving easily degenerates into ``constructing an objective function and then approximating it.'' Separating solving means that solving activities return to the original purpose of ``solving the problem,'' rather than returning to ``minimizing some external metric.''

Evaluation metrics and experience enter after the Solving Element. After construction and solving are completed, residuals, errors, or other metrics can still be calculated, but their role is verification, not guidance. If the structure is sufficient, the Solving Element executes the mapping in one step; if the structure is partially sufficient, the Solving Element executes a sequence of structural projections; if the structure is lacking, the Construction Element has no basis, and the Solving Element degenerates into Objective-Driven iteration. At this point, the Four-Element Theory degenerates into the Objective-Driven framework, rather than negating it.

The same applies to experience: it no longer participates in the construction and execution of the solving path, but only serves as a reference in the result verification stage. This positional adjustment may seem subtle, but it fundamentally cuts off the channel of distorting problem attributes to cater to metrics.

\section{Logical Constraints and Applicable Boundaries of the Four-Element Theory}

\subsection{Scientific Legitimacy Based on Falsificationism}

From the perspective of Popper's falsificationism, a comparison can be made between the working mechanisms of the two paradigms.

In the Objective-Driven paradigm, the evaluative objective (loss function, residual metrics) is placed at the core of the solving process. When systematic contradictions arise between the solution and real observations, the internal debugging path of the paradigm often prioritizes adjusting regularization, hyperparameters, or optimization strategies, rather than questioning the intrinsic mathematical-physical structure of the problem itself. This creates a practical risk: errors at the ontological structure level are masked as algorithmic tuning issues, making it difficult to refute structural hypotheses—though in principle, falsifiability is not entirely absent.

In contrast, the Structure-Driven paradigm adopts ``theory in the foreground, objective in the background'': the rigid invariants extracted after deconstruction are treated as conjectural hypotheses about the objective problem, upon which the solving mapping is constructed; evaluation metrics withdraw from the solution generation process and serve solely as \textit{a posteriori} testing tools. When the solution significantly conflicts with observations, the contradiction is first directed at this set of rigid invariant conjectures. Of course, according to the Duhem-Quine thesis, the mismatch of results can stem from multiple links such as deconstruction, construction, numerical implementation, and observational noise, and cannot be simply and solely attributed to the structural hypothesis. However, this process pushes the objective structural hypothesis to the forefront of testing, reducing the risk of structural errors being masked by parameter tuning, thereby better facilitating the scientific cognitive process of conjecture and refutation in Popper's sense.

\subsection{The Four Elements Cannot Be Reversed}
There are strict sequential constraints among the four stages. The Deconstruction Element must come first: structure identification occurs between basic units, and if the problem has not been reduced to independent units, the identified relationships will inevitably be mixed with coupling. The Structure Element must precede the Construction Element: the legitimacy of construction comes entirely from rigid invariants, and if structure identification is skipped and construction is done directly, construction degenerates into empirical guessing. The Construction Element must precede the Solving Element: the legitimacy of solving comes from the constructed mapping, and if there is no constructed mapping and solving is done directly, solving degenerates into search.

Thus, a criterion can be given: any solving method, as long as each step of its path can be traced back to some identified rigid invariant, and the final solving is executed by the mapping given by the Construction Element, it holds in the sense of Structure-Driven methodology; conversely, if there are untraceable steps in its path, or if solving is still guided by the gradient of the objective function, then no matter how sophisticated its form, it has not departed from the scope of Objective-Driven.

To determine whether a solving method belongs to the Structure-Driven Four-Element Theory, the following questions must be answered sequentially:
\begin{enumerate}
    \item Has a deconstruction suspending preset objectives been carried out to strip unit coupling as much as possible?
    \item Can each step of the solving path be traced back to a rigid invariant that has passed perturbation testing?
    \item Is the solving mapping constructed from the rigid structure, rather than dominated by the gradient of the objective function?
    \item Are evaluation metrics used only for posterior verification, without guiding the generation of the solving path?
\end{enumerate}
If all are satisfied: Structure-Driven; if partially satisfied: partially Structure-Driven; if mainly relying on objective gradient guidance: Objective-Driven.

\subsection{Applicable Boundaries}
The effectiveness of the Four-Element Theory depends on whether the problem has identifiable rigid structures. This condition constitutes the applicable boundary of the methodology.

When the structure is sufficient, solving can be completed completely directly, and both efficiency and interpretability are improved. When the structure is partially present, solving manifests as a sequence of structural projections, each step having a structural basis. When the structure is lacking, the methodology degenerates: at this time, the problem has no rigid invariants available for use, the Construction Element has no basis, and the Solving Element degenerates into Objective-Driven iteration.

This degeneration is not a defect of the methodology, but its honesty. Any methodology has its applicable boundaries, and acknowledging boundaries is closer to science than claiming omnipotence. More importantly, the degeneration relationship reveals the relationship between the two methodologies: Structure-Driven and Objective-Driven are not substitutes, but complements. The former provides a direct path for problems with sufficient structure, and the latter retains general solving capability for problems lacking structure.

However, behind this complementarity hides a profound cognitive cost: the reason Objective-Driven has long been popular is its lower cognitive threshold---one only needs to set an objective function, and the computational process can be handed over to the machine for iteration; while Structure-Driven requires the solver to possess strong domain intuition and mathematical literacy. Therefore, Structure-Driven is not superior in all scenarios; it is prepared for those solvers who are willing and able to pay the ``cognitive cost.''

\subsection{The Position of Experience}
The treatment of experience by the Four-Element Theory requires special explanation, because it is easily misread as a negation of experience.

This methodology does not negate the value of experience. Experience is precious knowledge accumulated by humanity in long-term solving practice, and it plays an irreplaceable role in problem formulation, structural conjecture, and result rationality judgment. In fact, structure identification itself often requires profound domain experience as guidance; the reason experts can see the structure is precisely the result of long-term internalization of experience.

What this methodology changes is the position of experience. The significance of this adjustment is reflected in three aspects:
\begin{itemize}[leftmargin=2em]
    \item \textbf{Shift in the source of legitimacy}: The source of legitimacy of the solving path is shifted from ``experience says so'' to ``structure says so.'' The former is non-reproducible because experience varies from person to person; the latter is reproducible within a given formalized framework because the structure is invariant for the same problem.
    \item \textbf{Relegation of the role}: Experience does not participate in the construction and execution of the solving path, but only after the solution is generated, is used to judge whether the result is reasonable and whether structural information needs to be supplemented.
    \item \textbf{Cutting off the core root cause}: This positional adjustment may seem subtle, but it fundamentally cuts off the channel of distorting problem attributes to cater to external metrics.
\end{itemize}

\subsection{Complementarity}
In scenarios where the structure is sufficient, the Structure-Driven path provides a more reliable solution at a lower cost, and the Objective-Driven path is not worse, but unnecessary. In scenarios where the structure is lacking, the Objective-Driven path is the only feasible general approach, and the Construction Element of the Structure-Driven path has no basis. In scenarios where the structure is partially sufficient, the two can collaborate: Structure-Driven handles the structurable part, and Objective-Driven handles the remaining degrees of freedom.

From the perspective of methodological evolution, Structure-Driven does not aim to replace Objective-Driven, but to fill the gap left by it. The reason the Objective-Driven methodology has been dominant across two centuries is precisely its unparalleled simplicity---one only needs to set an objective function, hand the solving process to the machine for iteration, and it yields a solution. However, its limitation is precisely the cost of ``sacrificing structure for generality.'' The Structure-Driven methodology sacrifices a part of generality in exchange for directness and interpretability in scenarios with sufficient structure. The two happen to be complementary.

\section{Positioning in the Methodological Spectrum}

The following discussion is not a deduction from philosophical concepts, but a methodological positioning extracted from the above practices.

\subsection{Relationship with Structuralism}
The idea of Structure-Driven has long roots in the humanities and social sciences. The structuralist movement of the 20th century argued that the meaning of phenomena does not depend on individual elements, but on the relational structure between elements; the research task is therefore not to describe elements, but to reveal structures. This is highly consistent with the core idea of the Structure Element: both believe that the essence of a problem lies in relationships rather than elements.

However, there are two differences between the Four-Element Theory and classical structuralism. First, the structures of classical structuralism are mostly speculative, and their existence is difficult to verify independently; the Four-Element Theory requires structures to be verifiable rigid invariants. Second, classical structuralism usually stops at structure revelation and does not involve how to generate results based on structure; the Four-Element Theory connects structure identification, construction generation, and solving execution into a complete chain.

\subsection{Relationship with Reductionism and Systems Theory}
As far as the Deconstruction Element is concerned, it has a superficial similarity with reductionism: both split complex objects into more basic units. But the goal of reductionism is to explain the properties of the whole with the properties of the units, and its implicit premise is that the whole can be completely determined by the parts; the goal of the Deconstruction Element is not explanation, but to create conditions for structure identification. The Four-Element Theory acknowledges that the whole has an irreducible relational structure, which is precisely its distinction from reductionism.

As far as the Structure Element is concerned, it has something in common with systems theory: both emphasize relationships and wholeness. But systems theory usually focuses on the general properties of systems, such as openness, hierarchy, and self-organization; the Four-Element Theory requires the extraction of specific, verifiable invariants, shifting the focus from general properties to specific constraints.

\subsection{Relationship with TRIZ, PINN, and Bayesian Priors}
Theory of Inventive Problem Solving (TRIZ) and other heuristic methods attempt to transcend empirical trial and error, forming reusable innovation principles by inducing solution patterns from a large number of patents. Its motivation is similar to the Four-Element Theory: both attempt to liberate solving from blind trial and error. However, TRIZ relies on cross-case induction, and its knowledge comes from the statistical regularities of existing solutions; the Four-Element Theory relies on the deduction of the problem's own structure, and its basis comes from constraints directly derived from the problem definition, not relying on a case base.

Physics-Informed Neural Networks (PINN) attempt to embed physical laws into the learning process \cite{ref12}, and they also intersect with the Four-Element Theory: both attempt to let structural constraints participate in solving. The difference is that the former treats structure as a soft constraint in the loss function, and the solving path is still Objective-Driven; the latter treats structure as the construction basis of the solving path, and structure enters the path itself rather than the objective function.

Bayesian Maximum A Posteriori estimation is a typical representative of Objective-Driven: it encodes structural information through prior distributions, combines it with the likelihood term into an objective function, and solves through optimization. The prior distribution is only a ``soft'' statistical expression of structural information, rather than an explicit utilization of the problem's intrinsic structure. The difference of the Structure-Driven methodology is that it does not ``stuff'' structure into the objective function, but ``extracts'' structure to construct inverse operators or solving mappings.

\subsection{Methodological Positioning}
The Four-Element Theory is closer to a research program in Lakatos's sense \cite{ref14}, rather than an incommensurable revolution in Kuhn's sense \cite{ref13}: it has a hard core (Structure-Driven, rigid invariants, construction-solving), a protective belt (specific structure identification methods, mathematical and physical structure paths, projection construction), and heuristics (starting from structure, rather than from objectives). Therefore, it coexists with Objective-Driven in both complementarity and competition, rather than complete replacement.

\section{Common Logic in Different Practices}

The following four cases belong to mathematics, physics, artificial intelligence, and geophysics, with vastly different problem types, object natures, and methodological traditions. If the Four-Element Theory exhibits completely consistent execution logic in them, its generality gains empirical support. Each case is presented in the same structure: what the problem is, how Objective-Driven does it, what Structure-Driven identified, how this structure makes solving direct, and what methodological proposition this case illustrates.

\subsection{Number Theory: From Asymptotic Approximation to Structural Proof}
\textbf{What the problem is.} In number theory, there is a class of propositions that require proving that certain prime pairs necessarily exist in even domains. The solving goal of such problems is to construct a rigorous proof, rather than to obtain a numerical approximation.

\textbf{How Objective-Driven does it.} The traditional path takes ``highest accuracy of asymptotic estimation of prime distribution'' as the objective, continuously increasing the order and accuracy of estimation, attempting to approximate the conclusion through more refined asymptotic analysis. The difficulty of this path lies in: no matter how refined the asymptotic estimation is, it only gives a probabilistic judgment of ``should be so,'' and cannot give a proof of necessity. There is a type gap between the accuracy metric and the proof goal: the former is analytical and approximate, while the latter is deductive and rigorous. Continuously raising accuracy cannot bridge this gap.

\textbf{What Structure-Driven identified.} The Structure-Driven path first deconstructs: reducing the global even domain to independent subdomains, stripping away the global averaging effect brought by asymptotic estimation in analytic number theory. Then it extracts rigid invariants: the multiplicative properties of completely multiplicative functions, and the commutative law of finite field multiplicative groups. Such properties do not change with specific numerical values and are constraints directly derived from the proposition definition.

\textbf{How the structure makes solving direct.} The Construction Element constructs a proof mapping accordingly; the Solving Element executes the proof mapping, deriving the conclusion step by step from structural constraints without any asymptotic approximation. At this time, proof logical self-consistency exists as a posteriori verification, rather than as a construction basis.

\textbf{Methodological proposition.} The key to this case is not number theory techniques, but the anchor of the solving path. The Objective-Driven path is anchored to ``highest asymptotic estimation accuracy,'' with the anchor on the solver's side; the Structure-Driven path is anchored to ``complete multiplicativity'' and ``commutative law,'' with the anchor on the problem's side. Once the anchor moves, the proof path changes from ``approximation'' to ``construction.''

The formal proof of the weakened form of the Liouville function Goldbach conjecture in 2026, and Wang Hong and Zahl's proof of the three-dimensional Kakeya conjecture, both embody this logic. It should be noted that the former has not yet undergone peer review and is only used here as an illustrative trend.

\subsection{Physics: From Error Minimization to Structure Preservation}
\textbf{What the problem is.} Long-term numerical simulation of high-dimensional Hamiltonian systems is a typical forward engineering problem. The core feature of such systems is that their phase space has a symplectic structure: the evolution process preserves phase space volume, and energy is conserved over long time scales.

\textbf{How Objective-Driven does it.} The traditional methodology takes minimizing single-step simulation error as the objective, continuously debugging the order and coefficients of numerical schemes, attempting to make the truncation error of each step as small as possible. However, this path has a structural blind spot: it excludes the symplectic structure from the evaluation system, because the symplectic structure cannot be characterized by single-step error metrics. Error metrics measure local deviation at a single step, while symplectic structure describes global geometric constraints; the two do not correspond in type.

\textbf{What Structure-Driven identified.} The Structure-Driven path first deconstructs: reducing the complete dynamical system to independent motion subunits, isolating non-structural coupling caused by non-conservative forces or high-frequency oscillations. Then it extracts rigid invariants: the symplectic symmetry of system evolution. It is invariant under any canonical transformation, and does not change with step size, scheme, or initial conditions.

\textbf{How the structure makes solving direct.} The Construction Element constructs structure-preserving numerical schemes accordingly (such as symplectic geometric algorithms \cite{ref11}); the Solving Element executes the scheme, making each step of evolution proceed under the constraints of the symplectic structure. At this time, there is no need to guide parameter tuning with error metrics; the form of the scheme is uniquely determined by the structure. The elimination of long-term energy drift is not tuned out, but guaranteed by the structure.

\textbf{Methodological proposition.} This case shows: the failure of the Objective-Driven path is not because it does not calculate finely enough, but because its evaluation metrics and the problem's rigid structure are not of the same type. Single-step error is a local quantity, and symplectic structure is a global geometric constraint; the former can be arbitrarily refined but never preserves structure. The reason the Structure-Driven path is effective is that it moves the anchor from ``minimum error'' to the invariant of ``symplectic symmetry.''

Noether's theorem has long revealed the profound connection between symmetry and conservation laws. Physicists have always been seeking the symmetry structures of problems, and then using them to simplify solving. This is precisely the historical paradigm of the Structure Element.

\subsection{Deep Learning: From Noise Fitting to Signal Fitting, from Gradient Descent to Structure Propagation}
\textbf{What the problem is.} Deep learning is the most typical Objective-Driven field of contemporary times. Whether it is image generation, signal denoising, or network training, the mainstream approach is to first define a loss function, then calculate gradients through backpropagation, and iteratively update parameters along the gradient direction. The form of the objective function is changing, the form of the optimizer is changing, but the core logic of ``constructing an objective function and minimizing it'' has not.

\textbf{How Objective-Driven does it.} In seismic signal denoising, traditional diffusion models restore signals by gradually adding noise and then gradually denoising, and their training objective is to fit the noise distribution. In neural network training, backpropagation calculates gradients through the chain rule, and then updates parameters along the gradient direction. The common point of the two is: the solving path is determined by the gradient of the objective function, not by the problem's own structure.

\textbf{What Structure-Driven identified.} The Structure-Driven path unfolds at two levels.

At the first level, in diffusion model denoising, when guiding a doctoral student to construct a diffusion model, the author proposed changing the noise fitting objective of the general diffusion model to a signal fitting objective based on the characteristic that seismic wave signal variation features are relatively simple, and led the naming of the method \cite{ref10}. The rigid structure identified here is: the structure of the seismic signal can be characterized by a low-dimensional representation of the signal itself, without relying on fine fitting of the noise distribution.

At the second level, in neural network training, the author proposed PSDBP \cite{ref9}, attempting to reconstruct backpropagation from the ``mechanical execution of gradient search'' to a ``propagation path constructed according to network structure.'' Here, the Deconstruction Element splits the network's forward propagation into independent hierarchical feature transformations, stripping the coupling between parameters of different layers; then it extracts the rigid structure: the forward propagation of the network has a decomposable causal chain and structural invariants. \textbf{This is a typical ``quasi-physical process'': we analogize the layer-by-layer transmission of information to the layer-by-layer propagation of physical wavefields.} Backpropagation can therefore be constructed according to these structures, without relying solely on the layer-by-layer gradient calculation of the chain rule.

\textbf{How the structure makes solving direct.} In diffusion model denoising, the Construction Element constructs a signal fitting mapping accordingly, and the Solving Element executes the mapping to directly restore the signal. This adjustment may seem subtle, but it moves the anchor of the solving path from ``noise distribution'' to ``signal structure.'' In backpropagation, the Construction Element constructs a structure-driven propagation mapping accordingly, and the Solving Element executes the mapping to directly complete parameter updates. The common logic of both is: to move the anchor from ``objective function'' to ``problem structure.''

\textbf{Methodological proposition.} This case shows: even in the field of deep learning, which is absolutely dominated by Objective-Driven, Structure-Driven can also emerge. It does not aim to abandon the objective function, but to change the anchor---from ``what to fit'' to ``what is the structure of the problem itself.'' This is completely consistent with the logic in number theory, physics, and geophysics. More importantly, PSDBP shows: even the core training algorithm of deep learning---backpropagation---can be reinterpreted from a Structure-Driven perspective. This shows that Structure-Driven is not a methodology external to artificial intelligence, but can enter its core mechanism.

\subsection{Geophysics: From Gradient Search to Structure Projection}
\textbf{What the problem is.} Seismic wave inversion imaging of strongly scattering complex media and 3D potential field inversion are typical inverse problems. The goal is to infer the elastic parameter or density distribution of the subsurface from wavefield or gravity field records observed at the surface.

\textbf{How Objective-Driven does it.} The traditional methodology takes minimizing the residual of observed data fitting as the objective, continuously updating the model through gradient iteration. Gradients are calculated by the adjoint state method, and iterations are driven by steepest descent, conjugate gradient, or quasi-Newton methods. The threefold limitations of this path are extremely typical here: each iteration requires a complete wavefield or potential field forward modeling, resulting in high computational cost; the gradient direction only indicates how to change the model to reduce residuals, and does not directly correspond to the physical causes of the data, limiting interpretability; the cross-correlation operation in gradient calculation mixes the contributions of different parameters, creating mathematical coupling, which requires the Hessian matrix or approximation methods to resolve, and the latter is difficult to implement in high-dimensional situations. More seriously, under strong scattering conditions, the objective function is highly non-convex, gradient iteration is extremely prone to falling into local minima, and the final imaging result deviates from the true subsurface structure.

\textbf{What Structure-Driven identified.} Under the Structure-Driven path, the two inversions exhibit extremely clear mathematical and physical structures, respectively. As a specific form of the Four-Element Theory in inverse problems, the two paths of SDI are reflected here:

For \textbf{true 3D gravity inversion} (MSDI): The Deconstruction Element splits the spatial data through \textbf{spectral-domain unitary transformation}, obtaining the basic components of the problem; it then identifies that the forward operator has a translation-invariant convolution property, recognizing its underlying \textbf{quasi-SVD structure}. This quasi-SVD structure is the rigid invariant of the problem.

For \textbf{strongly scattering seismic wave imaging} (PSDI/PSDWII): The Deconstruction Element reduces the subsurface medium space to independent scattering subunits, utilizing wavefield separation or time windows to strip highly coupled terms such as direct waves, scattered waves, and multiples under strong scattering conditions, and extracting independent incident wave and scattered wave units. Then it extracts physical rigid structures: the causal chain of wave propagation (incident wave---virtual source---secondary wave) and the decomposability of wavefield propagation.

\textbf{How the structure makes solving direct.} For gravity inversion: The Construction Element utilizes the quasi-SVD structure to directly construct a \textbf{pseudo-inverse mapping} in the spectral domain; the Solving Element executes this pseudo-inverse in one step, completing true 3D density imaging. For seismic wave imaging: The Construction Element directly constructs an imaging mapping through \textbf{virtual source projection} and deconvolution; the Solving Element executes the imaging mapping, generating the result. The local minimum problem is avoided not because a better optimization algorithm is used, but because there is no search link in the solving path at all.

\textbf{Methodological proposition.} This case most typically demonstrates the threefold limitations of Objective-Driven, and also demonstrates the two implementation paths of Structure-Driven in inverse problems. \textbf{Whether it is the ``quasi-SVD structure'' utilized by gravity inversion or the ``virtual source propagation physical structure'' utilized by seismic imaging, the essence is to recognize the skeleton of the problem through deconstruction, and then directly construct the solution.} The way the Structure-Driven path avoids them is not ``changing to a better optimization algorithm,'' but ``making the solving path have no search link at all.'' This is not a technical improvement, but a methodological shift.

\subsection{Consistency of the Four Cases}
The object natures, problem types, and methodological traditions of the four scenarios are vastly different, but under the framework of the Four-Element Theory, they exhibit completely consistent execution logic: deconstruction reduction, structure extraction, construction mapping, and solving execution.

More noteworthy is the consistency of the failure mechanism. In all four scenarios, the difficulties of the Objective-Driven methodology stem from the same structural blind spot: in the dynamical system scenario, error metrics obscure the symplectic structure; in the number theory scenario, accuracy metrics obscure number-theoretic invariants; in the deep learning scenario, noise fitting objectives obscure signal structure and network structure; in the inversion scenario, residual metrics obscure operator spectral structure and wavefield causal chain. And in all four scenarios, the elimination of difficulties follows the same path: not making metrics more refined, but identifying the structure.

\section{Systematic Comparison between Objective-Driven and Structure-Driven Paradigms}

\begin{table}[htbp]
\centering
\renewcommand{\arraystretch}{1.3}
\begin{tabular}{>{\raggedright\arraybackslash}p{0.22\textwidth} >{\raggedright\arraybackslash}p{0.35\textwidth} >{\raggedright\arraybackslash}p{0.35\textwidth}}
\toprule
\textbf{Comparison Dimension} & \textbf{Objective-Driven Methodology} & \textbf{Structure-Driven Four-Element Theory} \\
\midrule
Core Cognitive Anchor & Experience-defined preset objective evaluation system & Native rigid invariants inherent to the problem itself \\
\addlinespace
Core Solving Logic & Skips structure identification, iteratively tests and errors with the objective as the core & First completes structure identification, constructs forward mapping with rigid structure as the core, and executes the mapping to directly generate the solution \\
\addlinespace
Role of Experience & Pre-defined core constraints defining the solving boundary & A posteriori verification auxiliary reference after solving is completed \\
\addlinespace
Core Logical Root Cause & Ignores the inherent rigid structure of the problem from the root & Takes identifying the inherent rigid structure of the problem as the core starting point of solving \\
\addlinespace
Local Optimum Risk & Exists; complex nonlinear problems under objective constraints easily fall into local optima & No objective function; under the premise of no false invariants, since there is no objective-gradient search link, local optima issues of the Objective-Driven framework are avoided \\
\addlinespace
Result Reproducibility & Relatively low; results highly depend on initial models, hyperparameters, and solver's experience & When formalization and implementation conditions are consistent, the solving results for the same problem are completely consistent \\
\addlinespace
Cross-Domain Transferability & Limited, highly dependent on manually set objective functions and empirical tuning for specific scenarios (but widely used as a general means due to its high implementation simplicity) & Broad; solving logic constructed based on objective rigid invariants can be seamlessly transferred to different disciplines \\
\addlinespace
Resource Consumption Characteristic & High; massive computing power and manpower are used for empirical repeated trial and error & Low; only one mapping construction and execution is needed to obtain the final result \\
\bottomrule
\end{tabular}
\caption{Systematic comparison between the Objective-Driven methodology and the Structure-Driven Four-Element Theory.}
\label{tab:comparison}
\end{table}

The differences in the eight dimensions (see \autoref{tab:comparison}) can be attributed to one point: the position of the anchor. The anchor of the Objective-Driven methodology is on the solver's side, defined by experience; the anchor of the Structure-Driven methodology is on the problem's side, given by structure. The movement of the anchor position consequently triggers a series of differences in solving logic, role of experience, risk characteristics, reproducibility, transferability, and resource consumption.

\section{Methodological Significance, Boundaries, and Prospects}

\subsection{Dual Improvement in Efficiency and Interpretability}
The value of the Four-Element Theory is first reflected in the simultaneous improvement of efficiency and interpretability. In the Objective-Driven methodology, the two often conflict: improving interpretability usually means introducing more constraints, thereby sacrificing efficiency; improving efficiency usually means simplifying models, thereby sacrificing interpretability.

The reason the Structure-Driven path can improve both simultaneously is that it changes the source of both. Efficiency comes from direct construction and execution: solving is no longer trial-and-error iteration, but the mapping is directly given by the structure and executed. Interpretability comes from the structure itself: each step of the solving path can be traced back to an invariant with clear mathematical or physical meaning. When efficiency and interpretability originate from the same structure, the tension between them disappears.

From the methodological level, the emergence of Structure-Driven reflects a cognitive shift: from ``how to search'' to ``how to identify.'' For the past half century, the mainstream approach to inverse problem solving has been to construct error functions and iteratively search; the Structure-Driven approach is to utilize the reversible structure of the problem itself to directly construct a stable approximation of the inverse operator. Its core tools shift from differentiation to algebra and causal analysis. The significance of this shift lies not only in efficiency improvement but also in cognitive deepening: it makes us ask ``what is the essential structure of the problem,'' rather than ``how to search faster.''

\subsection{Theoretical Boundaries}
The boundary of the methodology has been given in Section 4: structure identifiability is the premise of effectiveness. Three additional qualifications are added here.

First, the existence of structure does not equal the identifiability of structure. The structures of some problems exist in principle, but are difficult to express explicitly under current mathematical tools. At this time, the methodology is theoretically applicable but practically limited. Second, the sufficiency of structure determines the degree of directness. Completely sufficient structures allow one-step construction and solving; partially sufficient structures require a sequence of structural projections. Therefore, directness is a matter of degree, not all-or-nothing. Third, structure identification itself may become a new bottleneck. When the problem dimension is extremely high and the unit relationships are extremely complex, the cost of structure identification may exceed the cost of iterative solving.

In addition, attention must be paid to the correctness of structure identification. If the so-called invariants identified are actually not inherent constraints of the problem, then the Construction Element will generate a mapping based on a false premise, and the Solving Element will execute the wrong mapping. The consequences are more serious than the local optima of the Objective-Driven methodology---the latter is at worst not reaching the optimum, while the former may give a systematically wrong solution. Therefore, structure identification must be premised on strict invariance testing, and cannot rely solely on analogy or intuition.

\subsection{Cross-Disciplinary Significance}
The significance of the Four-Element Theory is not limited to the methodological level. It reveals a common cross-disciplinary logic: engineering optimization, mathematical proof, physical simulation, and artificial intelligence training---these seemingly unrelated solving activities---share the same execution structure of deconstruction, structure, construction, and solving at the methodological level. This discovery helps break down methodological barriers between disciplines, enabling researchers in different fields to share a set of underlying solving language.

From a longer time scale, this methodological shift also has certain epistemological significance: it shifts the focus of solving activities from how to approach the objective faster to how to more completely understand the problem and grasp the structure. The former is a typical form of instrumental rationality, while the latter is closer to the exploration of the essence of things. This shift in the focus of solving activities may precisely be the methodological preparation humanity needs when facing increasingly complex systems.

\subsection{Future Work}
The verification boundary of this paper is strictly limited to the natural sciences and engineering fields, and the number of selected cases is limited, which is not yet sufficient to constitute statistical empirical evidence. In addition, due to the open nature of structure identification, this paper only provides a preliminary path induction for the general method of structure identification, and has not yet formed a complete operational specification, which is precisely the key to the implementation of the methodology.

Future work will proceed in three directions: first, expand cases to cover more problem types and disciplinary fields; second, study generalized methods of structure identification and explore whether there are cross-domain universal structure discovery strategies; third, examine projection construction and solving execution in situations where structure is incomplete, and clarify their stability conditions and error propagation laws.

\section{Conclusion}

Starting from concrete practices in natural science and engineering applications, this paper proposes and demonstrates the Four-Element Theory of Structure-Driven Methodology. The main conclusions are as follows.

First, the inherent root cause of the Objective-Driven methodology is not iteration and trial and error itself, but the pre-setting of the evaluative objective. Starting from experience-constrained objective evaluation inevitably leads to the systematic obscuring of the problem's internal structure; once the structural blind spot is formed, iterative trial and error becomes the only feasible path, and the threefold dilemma of local optima, non-reproducible results, and distorted native attributes follows.

Second, the Four-Element Theory reconstructs the solving process through four irreversible stages: deconstruction, structure, construction, and solving. The Deconstruction Element suspends \textit{a priori} presuppositions and reduces the problem to mutually independent basic units; the Structure Element extracts rigid invariants inherent to the problem; the Construction Element constructs equivalent mappings with the invariants as anchors; the Solving Element executes the mapping to directly generate the final solution. The fundamental basis of the four-element division is: deconstruction, structure, and construction all prepare for solving, and the Solving Element is the implementation stage of problem resolution. Evaluation metrics and experience are relegated to the position of \textit{a posteriori} verification.

Third, the applicable condition of this methodology is that the problem has identifiable rigid structures. When structure is sufficient, it can be directly constructed and solved; when structure is partially present, it manifests as a sequence of structural projections; when structure is lacking, it degenerates into the Objective-Driven path. The two methodologies are complementary rather than substitutes.

Fourth, the methodology exhibits consistent execution logic in different scenarios such as structured proofs of number theory propositions, long-term evolution simulation of Hamiltonian systems, deep learning signal processing and structure-driven backpropagation, seismic wave inversion imaging, and gravity inversion, and its generality gains cross-domain support. \textbf{In particular, as a specific form of the Four-Element Theory in inverse problems, SDI, along with its derived spectral-domain unitary diagonalized pseudo-inverse (based on quasi-SVD structure), true 3D gravity inversion, PSDWII (a true physical process), and PSDBP (a quasi-physical process), constitute the most representative practical foundation of the Four-Element Theory in mathematical inversion and deep learning.}

Fifth, in the methodological spectrum, the Four-Element Theory takes the holistic insight of structuralism as an ideological resource, verifiability as the entry threshold for structure identification, construction generation as the design stage, and solving execution as the purpose stage, forming a solving path that neither relies on case induction nor gradient search.

Sixth, the argumentation method of this paper itself also embodies this methodology: it does not deduce methodology from philosophical concepts, but extracts common structures from different practices such as number theory, physics, deep learning, and geophysics. This shows that the Four-Element Theory is applicable not only to objects but also to its own exposition---it is itself a practice of ``deconstruction---structure---construction---solving.''

Finally, from the perspective of scientific epistemology, the Objective-Driven approach is a ``target-catering'' fitting logic, whereas the Structure-Driven approach is an ``accepting-testing'' falsification logic. By relegating the objective to the background, SD takes physical reality as the sole arbiter for testing structural hypotheses. This is not merely a shift in research methods, but a return to the true essence of science.

The Four-Element Theory is the generalization and systematization of the previous SDI paradigm at the level of general problem solving. Its ultimate pursuit is to make solving activities return as much as possible to the original logic of directly solving using the inherent laws of the system, rather than relying on blind iterative approximation. This is not only an improvement in efficiency but also a deepening of cognition.

What this paper proposes is a research program for analyzing solving activities, rather than a fully mature, ready-to-apply technical toolkit; the structure identification stage still has open challenges that await further research.


\end{document}